\documentclass[11pt]{article}
\usepackage{amsthm}
\usepackage{amsfonts}
\usepackage[T1]{fontenc}
\usepackage[latin1]{inputenc}
\usepackage{epsfig}

\usepackage[lflt]{floatflt}
\usepackage{graphicx}
\usepackage{amsmath}
\usepackage{amssymb}
\usepackage{fancyhdr}
\usepackage{color}
\usepackage[affil-sl]{authblk}
\usepackage{cancel}
\usepackage{subfigure}
\usepackage[english]{babel}
\usepackage{fancyhdr}
\usepackage{makeidx}
\usepackage{latexsym}

\newtheorem{theorem}{Theorem}

\newtheorem{lemma}{Lemma}
\newtheorem{corollary}{Corollary}
\newtheorem{remark}{Remark}

\makeindex \textwidth=17cm \textheight=25cm \oddsidemargin=-0.5cm
\newcommand{\bbR}{{\mathbb R}}

\newcommand{\bbC}{{\mathbb C}}

\newcommand{\al}{\alpha}
\newcommand{\la}{\lambda}
\newcommand{\p}{\partial}
\newcommand{\be}{\beta}
\newcommand{\G}{\Gamma}

\begin{document}
\begin{center}
\LARGE
\textbf{TWO EQUIVALENT STEFAN'S PROBLEMS FOR THE TIME FRACTIONAL DIFFUSION EQUATION.}
\end{center}
\begin{center}
\normalsize
 Sabrina Roscani$^1$ and Eduardo A. Santillan Marcus$^2$\\
\medskip
\small
$^{1,2}$Departamento de Matem\'{a}tica,  FCEIA, Universidad Nacional de Rosario, Pellegrini 250, Rosario, Argentina \\
\textcolor{blue}{greyero@fceia.unr.edu.ar, sabrina@fceia.unr.edu.ar, edus@fceia.unr.edu.ar} \\

$^2$ Departamento de Matem\'atica,
FCE, Universidad Austral, Paraguay 1950, Rosario, Argentina,\\
\textcolor{blue}{esantillan@austral.edu.ar}\\

$^1$ CONICET, Argentina.

\end{center}

\small

\noindent \textbf{Note: } This paper is now
published in Fract. Calc. Appl. Anal. Vol. 16, No 4 (2013), pp. 802-815, DOI: 10.2478/s13540-013-0050-7 , and is
available at http://link.springer.com/journal/13540.
\smallskip

\noindent \textbf{Abstract: }
Two Stefan's problems for the diffusion fractional equation  are
solved, where  the fractional derivative of order $ \al \in (0,1) $
is taken in the Caputo's sense. The first one has a constant
condition on $ x = 0 $ and the second presents a flux condition  $
T_x (0, t) = \frac {q} {t ^ {\al/2}} $. An equivalence between these
problems is proved and  the convergence to the classical solutions
is analyzed when  $ \al \nearrow $ 1 recovering the heat equation
with its respective Stefan's condition.
 \medskip


 \smallskip


 \section{Introduction}\label{sec:1}

\setcounter{section}{1}
\setcounter{equation}{0}\setcounter{theorem}{0}

 In 1695 L'H$\hat{\textrm{o}}$pital inquired of Leibnitz,
the father of the concept of the classical differentiation, what
meaning could be ascribed to the derivative of order $\frac{1}{2}$.
Leibnitz replied prophetically: ``[...] this is an apparent paradox
from which, one day, useful consequences will be drawn.''

Since 1819, mathematicians like Lacroix, Abel, Liouville, Riemann
and later Gr\"{u}nwald and Letnikov, have attempted to establish a
definition of derivative of fractional order.

Here we use the definition introduced by Caputo \cite{Caputo} in
1967, referred to as \textsl {fractional derivative in Caputo's
sense},  given by
$$\,^C_{a} D^{\alpha}f(t)=\frac{1}{\Gamma(n-\alpha)}\int^{t}_{a}(t-\tau)^{n-\alpha-1} f^{(n)}(\tau)d\tau, $$
where $\al >0$ is the order of differentiation, $n=\left\lceil \al
\right\rceil$ and  $f $ is a differentiable function up to order $n$
in $\left[a,b\right]$.

The one-dimensional heat equation has become the paradigm for the
all-embracing study of parabolic partial differential equations,
linear and nonlinear. Cannon \cite{Cannon} did a methodical
development of a variety of aspects of this paradigm. Of particular
interest are the discussions on the one-phase Stefan problem, one of
the simplest examples of a free-boundary-value problem for the heat
equation (see \cite{dat}). In mathematics and its applications,
particularly related to phase transitions in matter, a Stefan
problem is a particular kind of boundary value problem for a partial
differential equation, adapted to the case in which a phase boundary
can move with the time. The classical Stefan problem aims to
describe the temperature distribution in a homogeneous medium
undergoing a phase change, for example ice passing to water: this is
accomplished by solving the heat equation imposing the initial
temperature distribution on the whole medium, and a particular
boundary condition, the Stefan condition, on the evolving boundary
between its two phases. Note that in the one-dimensional case this
evolving boundary is an unknown curve: hence, the Stefan problems
are examples of free boundary problems. A large bibliography on free
and moving boundary problems for the heat-diffusion equation was
given in \cite{Tar2}.

\smallskip

In this paper, we study a one-phase Stefan problem with time
fractional diffusion equation, obtained from the standard diffusion
equation by replacing the first order time-derivative by a
fractional derivative of order $\alpha > 0 $ in the Caputo sense:
$$  \,^C_{a} D^{\alpha}u(x,t)=\lambda^2\dfrac{\partial^2 u }
{\partial x^2}(x,t), \quad  -\infty<x<\infty, \ t>0, \ 0<\al<1, $$
 and the Stefan condition $\frac{d s(t)}{dt}=k u_x(s(t),t), \, t>0,$ by the fractional Stefan condition
$$  \,^C_{a} D^{\alpha}s(t)= k u_x(s(t),t),\quad t>0. $$

 This equation has been recently treated by a number of
authors (see e.g. Gorenflo and Mainardi \cite{FM-Anal prop and apl
of the W func}, Liu and Xu \cite{Liu-Xu}, Kilbas \cite{Kilbas},
Podlubny \cite{Podlubny}) and, among the several applications that
have been studied, Mainardi \cite{FM-libro} studied the application
to the theory of linear viscoelasticity.

The solutions of this equation are expressed in terms  of two
special functions that play a very important role in the theory of
differentiation of arbitrary order: the Mittag-Leffler function (see
e.g. \cite{FM-libro}, \cite{Kilbas}, \cite{Podlubny})
$$E(z,\al,\be)=\sum^{\infty}_{n=0}\frac{z^n}{ \G\left( \al n+\be \right)}$$
and the Wright function \cite{Wright} 
$$W(z,\al,\be)=\sum^{\infty}_{n=0}\frac{z^n}{n! \G\left( \al n+\be \right)},$$
respectively. A particular case of the Wright function is the
 Mainardi function (see Podlubny \cite{Podlubny})
$$ M_\nu (z)= W(-z,-\nu,1-\nu)=\sum^{\infty}_{n=0}\frac{(-z)^n}{n! \G\left( -\nu n+ 1-\nu \right)}.  $$
 This function is a part of the fundamental solution for the time fractional diffusion equation
 studied in \cite{FM-The fundamental solution}
$$ G_\al(x,t)=\frac{1}{2\lambda t^{\al/2}}M_{\al/2}\left(\frac{x}{\lambda t^{\al/2}}\right). $$

\section{Solving the two Fractionary Stefan's problems}
\label{sec:2}

 Hereinafter, $\,^C_{0} D^{\alpha}=D^\al$. We consider the following
problems:
\begin{equation}{\label{St1}}
\left\{\begin{array}{lll}
          D^{\al} u(x,t)=\lambda^2\dfrac{\partial^2u }{\partial x^2}(x,t) &   0<x<s(t), \,  t>0, \,  0<\al<1 , \, \, \la>0\\
          u(0,t)=B   &  t>0 \quad  B>0 \ \mbox{constant} \\
          u(s(t),t)=C<B & t>0\\
          D^{\al}s(t)=-k u_x(s(t),t) & t>0, \quad  k>0\\
          s(0)=0                                   \end{array}\right.\end{equation}
and
\begin{equation}{\label{St2}}
\left\{\begin{array}{lll}
          D^{\al} u(x,t)=\lambda^2\dfrac{\partial^2u }{\partial x^2}(x,t) &  & 0<x<s(t), \, t>0, \, 0<\al<1, \, \, \la>0 \\
          u_x(0,t)=-\frac{q}{t^{\al/2}} \quad  \quad &   & t>0, \quad q>0 \\
          u(s(t),t)=C \quad  \quad & & t>0\\
          D^{\al}s(t)=-k u_x(s(t),t)\quad  \quad & & t>0
          \\
          s(0)=0                                   \end{array}\right. .
          \end{equation}

A pair $\{u,s\} $ is a solution of the problem $(\ref{St1})$ (or
$(\ref{St2})$) if:

\begin{enumerate}
    \item $u$ and $s$ satisfy (\ref{St1}) (or (\ref{St2})),
    \item $u_{xx} $ and  $D^{\al} u$ are continuous for $0<x< s(t)$, $0<t<T$,
    \item $u$ and $u_x$ are continuous for $0\leq x\leq s(t)$,
    $0<t<T$,
    \item $0\leq \underset{x,t\rightarrow 0^+}{\liminf}u(x,t)\leq \underset{x,t\rightarrow 0^+}{\limsup }
    u(x,t)<+\infty$,
    \item $s$ is continuously differentiable in $[0,T)$ and $\frac{\dot{s}(\tau)}{(t-\tau)^\al}$ $\in L^1(0,t)$
    $\forall t\in (0,T)$.
    \end{enumerate}

 Let us solve the problem (\ref{St1}). We show in Appendix that
 \begin{equation}\label{St1-i}u_1(x,t)=a_1+b_1\left[ 1-W\left(-\frac{x}{\la t^{\al/2}},-\frac{\al}{2},1\right)\right],
 \qquad  a_1,b_1 \, \text{ constant},
 \end{equation}
  is a solution for the time-fractional-diffusion equation.
\begin{equation}\label{St1-ii} u_1(0,t)=a_1+b_1\left[1-W\left(0,-\frac{\al}{2},1\right)\right]=B\Rightarrow a_1=B.
\end{equation}
\begin{equation}\label{St1-iii} u_1(s_1(t),t)=a_1+b_1\left[ 1-W\left(-\frac{s(t)}{\la t^{\al/2}},-\frac{\al}{2},1\right)\right]=C.
\end{equation}
 Note that (\ref{St1-iii}) must be verified for all $t>0$, so we will ask for  $s(t)$ to be proportional to $t^{\al/2}$,
 that is to say
\begin{equation}\label{s(t)}
s_1(t)=\la \xi t^{\al/2} \qquad  \text{for some } \xi>0,
\end{equation}
 and from (\ref{St1-iii}), (\ref{s(t)}) and Corollary
 \ref{1menosW>0},
\begin{equation}\label{St1-b}
C=B+b_1\left[ 1-W\left(-\xi,-\frac{\al}{2},1\right)\right]\Rightarrow b_1
 = \frac{C-B}{1-W\left(-\xi,-\frac{\al}{2},1\right)}\,.
 \end{equation}
Now we will obtain $\xi$ from the ``fractional Stefan condition''.
Taking into account that
$$ D^{\al}(t^{\be})=\frac{\G(\be+1)}{\G(1+\be-\al)}t^{\be-\al} \quad \text{if } \be>-1, $$
we have
\begin{equation}\label{St1-iii-a}
D^{\al}s_1(t)=D^{\al}(\la \xi t^{\al/2})=\la
\xi\frac{\G(\frac{\al}{2}+1)}{\G(1-\frac{\al}{2})}t^{-\al/2}.
\end{equation}
On the other hand,
 \begin{equation}\label{St1-iii-b}
 u_{1x}(s_1(t),t)=b_1 \frac{1}{\la t^{\al/2}}M_{\al/2}\left(\xi\right)=
 \frac{C-B}{1-W\left(-\xi,-\frac{\al}{2},1\right)}\frac{1}{\la t^{\al/2}}M_{\al/2}\left(\xi\right).
 \end{equation}
From $(\ref{St1-iii-a})$ and $(\ref{St1-iii-b})$
\begin{equation}\label{St1-eq para xi}
 \xi\left[1-W\left(-\xi,-\frac{\al}{2},1\right)\right]
 \frac{1}{M_{\al/2}\left(\xi\right)}=-\frac{k}{\la^2}
 \frac{\G(1-\frac{\al}{2})}{\G(1+\frac{\al}{2})}(C-B).
 \end{equation}
 Let us define
\begin{equation}\label{St1-H(xi)}
H(\xi)=\xi\left[1-W\left(-\xi,-\frac{\al}{2},1\right)\right]
 \frac{1}{M_{\al/2}\left(\xi\right)}.
 \end{equation}

The function $H$ has the following properties:

\begin{enumerate}
\item $H(0^+)=0,$
\item $H(+\infty)=+\infty,$
\item $H$ is continuous and monotonically increasing.
\end{enumerate}

Because of the asymptotic behavior of the Wright function (see \cite{FM-Anal prop and apl of the W func}),
 it is easy to check the Properties 1 and 2.

For Property 3, we observe from Corollary \ref{1menosW>0} that,
$1-W\left(-\xi,-\frac{\al}{2},1\right)$ is a positive and increasing
function in $\bbR^+$. And from Lemma \ref{M pos y decrec},
$\frac{1}{M_{\al/2}(\xi)}$ is a positive increasing function.

Observing that $-\frac{k}{\la^2}
 \frac{\G(1-\frac{\al}{2})}{\G(1+\frac{\al}{2})}(C-B)>0$, we can assure that there exists  a
 unique $\tilde{\xi}$ such that
\begin{equation}\label{H(xi)=algo} H(\tilde{\xi})=-\frac{k}{\la^2}
 \frac{\G(1-\frac{\al}{2})}{\G(1+\frac{\al}{2})}(C-B).
 \end{equation}
 So the  solution of problem (\ref{St1}) is given by
\begin{equation}{\label{SOL-St1}}
\left\{\begin{array}{lll}
           u_1(x,t)=B+\frac{C-B}{1-W\left(-\tilde{\xi},-\frac{\al}{2},1\right)}
           \left[1-W\left(-\frac{x}{\la t^{\al/2}},-\frac{\al}{2},1\right)\right]\\
s_1(t)=\la \tilde{\xi} t^{\al/2}, \\
          \text{where } \tilde{\xi} \text{ is the unique solution to the equation}  \\

          H(\xi)=-\frac{k}{\la^2}
 \frac{\G(1-\frac{\al}{2})}{\G(1+\frac{\al}{2})}(C-B).
            \end{array}\right.\end{equation}

 Now let us solve (\ref{St2}). Here we consider
 \begin{equation}\label{St2-i}u_2(x,t)=a_2+b_2\left[ 1-W\left(-\frac{x}{\la t^{\al/2}},-\frac{\al}{2},1\right)\right],
 \qquad  a_2,b_2 \, \text{constant}.
 \end{equation}
Then,
\begin{equation}\label{St2-b} u_{2x}(0,t)=\frac{b_2}{\la t^{\al/2}}M_{\al/2}(0)
=-\frac{q}{t^{\al/2}}\Rightarrow b_2=-q\la
\G\left(1-\frac{\al}{2}\right),
\end{equation}
\begin{equation}\label{St2-iii} u_2(s_2(t),t)=a_2+b_2\left[ 1-W\left(-\frac{s_2(t)}{\la t^{\al/2}},
-\frac{\al}{2},1\right)\right]=C.
\end{equation}
Note that (\ref{St2-iii}) must be verified for all $t>0$, so we will
ask for $s_2(t)$ to be proportional to $t^{\al/2}$, that is to say
\begin{equation}\label{s(t)-St2}
s_2(t)=\la \mu t^{\al/2}, \qquad  \text{for some } \mu
>0.
\end{equation}
From (\ref{St2-iii}) and (\ref{s(t)-St2}) we have
\begin{equation}\label{St2-a} a_2=C+q\la \G\left(1-\frac{\al}{2}\right)
\left[ 1-W\left(-\mu,-\frac{\al}{2},1\right)\right].
\end{equation}

Notice that
\begin{equation}\label{St2-iii-a}
D^{\al}s_2(t)=\la \mu
\frac{\G(\frac{\al}{2}+1)}{\G(1-\frac{\al}{2})}t^{-\al/2}
\end{equation}
and
 \begin{equation}\label{St2-iii-b}
 u_{2x}(s_2(t),t)=\frac{b_2}{\la
 t^{\al/2}}M_{\al/2}(\mu)=-\frac{q \G\left(1-\frac{\al}{2}\right)}{ t^{\al/2}} M_{\al/2}(\mu).
 \end{equation}
So, from $(\ref{St2-iii-a})$ and $(\ref{St2-iii-b})$,
$$ \la \mu \frac{\G(\frac{\al}{2}+1)}{\G(1-\frac{\al}{2})}t^{-\al/2}
=-k \frac{(-q)\G\left(1-\frac{\al}{2}\right)}{ t^{\al/2}}
 M_{\al/2}(\mu), $$
therefore
\begin{equation}\label{St2-eq para mu}
  \mu \frac{1}{ M_{\al/2}(\mu)}=\frac{kq}{\la}
 \frac{\G\left(1-\frac{\al}{2}\right)^2}{\G(\frac{\al}{2}+1)}\,.
 \end{equation}

Let us define
\begin{equation}\label{St2-J(xi)}
J(\mu)=\mu \frac{1}{ M_{\al/2}(\mu)}.
 \end{equation}
The function $J$ has the following properties:

\begin{enumerate}
\item $J(0^+)=0,$
\item $J(+\infty)=+\infty,$
\item J is continuous and monotonically increasing.
\end{enumerate}

Observing that $\frac{kq}{\la}
 \frac{\G\left(1-\frac{\al}{2}\right)^2}{\G(\frac{\al}{2}+1)}>0$,
  we can assure that there exists  a unique $\tilde{\mu}$ such that
$$  J(\tilde{\mu})=\frac{kq}{\la}
 \frac{\G\left(1-\frac{\al}{2}\right)^2}{\G(\frac{\al}{2}+1)}. $$
So the  solution of problem (\ref{St2}) is given by
\begin{equation}{\label{SOL-St2}}
\left\{\begin{array}{lll}
           \begin{array}{ll}
             u_2(x,t)= & C+q\la\G\left(1-\frac{\al}{2}\right)
             \left[ 1-W\left(-\tilde{\mu},-\frac{\al}{2},1\right)\right] \\
              & -q\la\G\left(1-\frac{\al}{2}\right)
           \left[ 1-W\left(-\frac{x}{\la
           t^{\al/2}},-\frac{\al}{2},1\right)\right],
           \end{array}
           \\
s_2(t)=\la \tilde{\mu} t^{\al/2}, \\
          \text{where } \tilde{\mu} \text{ is the unique solution to  equation}  \\

          J(\mu)=\frac{kq}{\la}
 \frac{\G\left(1-\frac{\al}{2}\right)^2}{\G(\frac{\al}{2}+1)}.
         \end{array}\right.\end{equation}

 Finally, our goal is to show the relationship between the two
diffusion fractional problems with temperature and flux conditions
at $x=0$, respectively, to obtain a similar result as the one given
by \cite{Tarzia}.

 \begin{theorem}\label{T2.1}
 Let us consider problems (\ref{St1}) and (\ref{St2}), where:
\smallskip

 (1) the constant $C$ is the same in both problems,

 (2) in problem (\ref{St1}):
  $B=C-q\la\G\left(1-\frac{\al}{2}\right)\left[ 1-W\left(-\tilde{\mu},-\frac{\al}{2},1\right)\right],$ where $\tilde{\mu}$ is the unique solution to $J(\mu)=\frac{kq}{\la}
 \frac{\G\left(1-\frac{\al}{2}\right)^2}{\G(\frac{\al}{2}+1)}$ and $J$ is defined by $(\ref{St2-J(xi)})$.

\smallskip

Then these problems are equivalent.
 \end{theorem} 

\proof
  Let us define the following function
$$ B(\xi)=C+q\la\G\left(1-\frac{\al}{2}\right)\left[ 1-W\left(-\xi,-\frac{\al}{2},1\right)\right].$$
Observe that $B(\xi)>C$, $\forall\,  \xi$ and $B(\tilde{\mu})=B$.

Now,
 $$H(\xi)=-\frac{k}{\la^2}
 \frac{\G(1-\frac{\al}{2})}{\G(1+\frac{\al}{2})}(C-B(\xi)) \Longleftrightarrow $$
$$ \xi\left[1-W\left(-\xi,-\frac{\al}{2},1\right)\right]
 \frac{1}{W\left(-\xi,-\frac{\al}{2},1-\frac{\al}{2}\right)}$$
$$=-\frac{k}{\la^2}\frac{\G(1-\frac{\al}{2})}{1+\frac{\al}{2}}(-q)a\G\left(1-\frac{\al}{2}\right)\left[
1-W\left(-\xi,-\frac{\al}{2},1\right)\right]$$
$$\Longleftrightarrow
\xi\left[1-W\left(-\xi,-\frac{\al}{2},1\right)\right]
=\frac{kq}{\la}\frac{\G(1-\frac{\al}{2})^2}{\G\left(1+\frac{\al}{2}\right)}\Longleftrightarrow$$
\begin{equation}\label{J teo eq}
\Longleftrightarrow
J(\xi)=\frac{kq}{\la}\frac{\G(1-\frac{\al}{2})^2}{\G\left(1+\frac{\al}{2}\right)}.
\end{equation}
 Then if $\tilde{\mu}$ is the unique solution of (\ref{J
teo eq}), we have
 $$H(\tilde{\mu})=\frac{k}{\la^2}
 \frac{\G(1-\frac{\al}{2})}{\G(1+\frac{\al}{2})}(C-B)$$
 $$=\frac{k}{\la^2}\frac{\G(1-\frac{\al}{2})}{\G(1+\frac{\al}{2})}q\la\G\left(1-\frac{\al}{2}\right)
 \left[ 1-W\left(-\tilde{\mu},-\frac{\al}{2},1\right)\right]>0. $$

Due to the uniqueness of solution of $(\ref{H(xi)=algo})$, we can
assure that $ \tilde{\mu}=\tilde{\xi} $, and therefore $ s_1=s_2 $.

 It is easy now to check that $u_1=u_2$. From (\ref{SOL-St1}) and
(\ref{SOL-St2})
$$  u_2(x,t)=C+q\la\G\left(1-\frac{\al}{2}\right)
\left[1-W\left(-\tilde{\mu},-\frac{\al}{2},1\right)\right]$$
$$-q\la\G\left(1-\frac{\al}{2}\right)\left[ 1-W\left(-\frac{x}{\la
t^{\al/2}},-\frac{\al}{2},1\right)\right] $$
$$=B-q\la\G\left(1-\frac{\al}{2}\right)\frac{1-W\left(-\tilde{\xi},-\frac{\al}{2},1\right)}
{1-W\left(-\tilde{\xi},-\frac{\al}{2},1\right)}\left[
1-W\left(-\frac{x}{\la t^{\al/2}},-\frac{\al}{2},1\right)\right]
 $$
$$=B+\frac{C-B}{1-W\left(-\tilde{\xi},-\frac{\al}{2},1\right)}
\left[1-W\left(-\frac{x}{\la t^{\al/2}},-\frac{\al}{2},1\right)\right]=u_1(x,t).
 $$
\endproof


\begin{remark}
Applying Theorem \ref{conv a la erf} from Appendix, to the given
solutions (\ref{SOL-St1}) and (\ref{SOL-St2}) we recover the
classical solutions:
$$\lim_{\al\nearrow 1}u_1(x,t)=\lim_{\al\nearrow 1}\left\{B+\frac{C-B}{1-W\left(-\tilde{\xi},-\frac{\al}{2},1\right)}
           \left[1-W\left(-\frac{x}{\la t^{\al/2}},-\frac{\al}{2},1\right)\right] \right\} $$
 $$ =B+\frac{C-B}{\mbox{erf\,}\left(\tilde{\xi}/2\right)}\, \mbox{erf\,}\left(-\frac{x}{2\la t^{1/2}}\right), $$
 $$\lim_{\al\nearrow 1}s_1(t)=\lim_{\al\nearrow 1}\la \tilde{\xi} t^{\al/2}=\la \tilde{\xi}\sqrt{t},
$$
where $ \tilde{\xi}$  is the unique solution to the equation
  $$\frac{\xi}{2} \mbox{erf\,}\left(\frac{\xi}{2}\right)e^{\xi^2/4}=
 \frac{k}{\la^2}\frac{(C-B)}{\sqrt{\pi}}. $$

\end{remark}

\section{Conclusions}\label{sec:3}

 We have studied the behavior of the two Wright functions in
$\bbR^+_0$: \break $1-W(-x,-\frac{\al}{2},1)$  and $M_{\al/2}(x)$,
and then we solved two fractional Stefan's problems for the time
fractional diffusion equation with its respective fractional
Stefan's conditions: the first one with a constant condition at
$x=0$, and the second one with a flux condition
$u_x(0,t)=-\frac{q}{t^{\al/2}}$. Finally, we proved the equivalence
between these two problems (for a suitable constant condition) and
we have analyzed the convergence when $\al\nearrow 1$, thus
recovering the classical solution to the heat equation and its
respective Stefan's condition.

\smallskip

\section{Appendix: Working with the Wright function}
\label{App}

 \setcounter{equation}{0}\setcounter{theorem}{0}

 Note that the Wright function
  $W(z,\al,\be)=\sum^{\infty}_{n=0}\displaystyle{\frac{z^n}{n! \G\left( \al n+\be \right)}}$
  is an entire function if $\Re(\al)>-1$.

  Taking $\al=-\frac{1}{2}$ and $\be=\frac{1}{2}$, we get
 $$  W\left(-z,-\frac{1}{2},\frac{1}{2}\right)=M_{1/2}(z)=\frac{1}{\sqrt{\pi}}e^{-z^2/4}.$$
 Due to the uniform convergence of the series,
 \begin{equation}\label{derivada de W} \frac{\p}{\p z} W(z,\al,\be) = W(z,\al,\al+\be).  \end{equation}
 Then, for $x\in \bbR^+_0$, and taking account that
 \begin{equation}\label{W(-inf,-al 2,1)} W(-\infty,-\frac{\al}{2},1)=0, \qquad  \text{ if } \
  \al\in (0,2),
  \end{equation}
 we have
 $$ W\left(-x,-\frac{1}{2},1\right)=  W\left(-x,-\frac{1}{2},1\right)-
 W\left(-\infty,-\frac{1}{2},1\right)$$
$$ = \int_{\infty}^x \left( \frac{\p}{\p x}
W\left(-\xi,-\frac{1}{2},1\right)\right)d\xi= \int_{\infty}^x
-W\left(-\xi,-\frac{1}{2},\frac{1}{2}\right)d\xi $$
$$ = \int_x^{\infty}
W\left(-\xi,-\frac{1}{2},\frac{1}{2}\right)d\xi=\int_x^{\infty}
\frac{1}{\sqrt{\pi}}e^{-\xi^2/4}d\xi$$
 $$ = \frac{2}{\sqrt{\pi}}\int_{x/2}^{\infty} \frac{1}{\sqrt{\pi}}e^{-\xi^2}d\xi
 = \mbox{erfc\,}\left(\frac{x}{2}\right). $$
 Consequently,
$$  W\left(-x,-\frac{1}{2},1\right)= \mbox{erfc\,}\left(\frac{x}{2}\right)$$
and
$$ 1-W\left(-x,-\frac{1}{2},1\right)= \mbox{erf\,}\left(\frac{x}{2}\right).$$


\begin{remark}
It is a fact that the Mainardi function $M_{\al/2}(z)$ is an entire function of $z$
 (see \cite{FM-Anal prop and apl of the W func}). So, any limit on the variable $z$ can be calculated
 by interchanging limit and sum. However that is not always true if the limit is taken in the parameter $\al$.

For example, the function $$ f_\al(z)=e^{-z/\al} $$ is an entire function on the variable $z$,
whose series representation is \break 
$\sum\limits_{n=0}^{\infty}\frac{(-z/\al)^n}{n!} $, and,  for every $z$ fixed,
$$ \lim\limits_{\al\searrow 0}e^{-z/\al}=0 $$
while
$$\lim\limits_{\al\searrow 0}  \frac{(-z/\al)^n}{n!}=\pm \infty, $$
and therefore the series diverges.
\end{remark} 


\begin{lemma}\label{conv M al/2 cuando al tiende a 1}
 If $x\in \bbR^+_0$ and $\al \in (0,1)$,
$$
\lim_{\al\nearrow
1}M_{\al/2}\left(x\right)=M_{1/2}(x)=\frac{e^{-\frac{x^2}{4}}}{\sqrt{\pi}}\,.
$$
\end{lemma} 

 \proof
 Let $\al$ be such that $\frac{1}{2}<\al< 1$.
 Writing the series as a sum of even and odd terms subseries, it will be seen that each  one of them
  is bounded by  a convergent series which does not depend on $\al$.
 For the even terms,
$$ \left|\sum_{k=0}^{\infty}
\frac{x^{2k}}{(2k)!\G(-\frac{\al}{2}2k+1-\frac{\al}{2})}\right|\leq
\sum_{k=0}^{\infty}\frac{|x|^{2k}}{(2k)!\left|\G(-\frac{\al}{2}2k+1-\frac{\al}{2})\right|}$$
$$=
\sum_{k=0}^{\infty}\frac{|x|^{2k}}{(2k)!\left|\G(1-\al(k+\frac{1}{2}))\right|}\,.
$$

Recall that for all $x\, \in  \bbR,$
$$\frac{1}{\G(x)\G(1-x)}=\frac{\sin (\pi x)}{\pi} ,$$
 and the Gamma function is increasing in $(\frac{3}{2},+\infty)$.
Then if $k\geq 3$,
$$ 0<\G(\al(k+\frac{1}{2}))\leq \G(k+1) \ \text{, \ \ and therefore } \ \
 \frac{1}{\G(k+1)}\leq \frac{1}{\G(\al(k+\frac{1}{2}))}\,.$$
 On the other hand, $(2k)!=(2k)...(k+1)\G(k+1)$. So,
$$  \frac{|x|^{2k}}{(2k)...(k+1)\G(k+1)\left|\G(1-\al(k+\frac{1}{2}))\right|}
$$
$$\leq \frac{|x|^{2k}}{(2k)...(k+1)\left|\G(\al(k+\frac{1}{2}))\G(1-\al(k+\frac{1}{2}))\right|}$$
$$=\left| \frac{|x|^{2k} \sin (\pi \al(k+1))}{(2k)...(k+1)\pi}\right|\leq \frac{|x|^{2k}k!}{\pi
(2k)!}, \quad \forall \, k\geq 3.  $$
Then
 \begin{equation}\label{conv unif serie term pares}
 \left|\sum_{k=0}^{\infty}
\frac{x^{2k}}{(2k)!\G(-\frac{\al}{2}2k+1-\frac{\al}{2})}\right|\leq
\left|\sum_{k=0}^{2}\frac{x^{2k}}{(2k)!\G(-\frac{\al}{2}2k+1-\frac{\al}{2})}\right|
+ \sum_{k=3}^{\infty} \frac{|x|^{2k}k!}{\pi (2k)!}\,.
 \end{equation}
It is easy to see that  this is an absolutely convergent series in $\bbC$.\\

Concerning the odd terms, reasoning in the same way, now with $k\geq
2$ we get
$$ \left|\sum_{k=0}^{\infty}
\frac{-x^{2k+1}}{(2k+1)!\G(-\frac{\al}{2}2(k+1)+1-\frac{\al}{2})}\right|
$$
$$\leq \left|\sum_{k=0}^{1}
\frac{-x^{2k+1}}{(2k+1)!\G(1-\al(k+1))}\right|+\sum_{k=2}^{\infty}\frac{|x|^{2k+1}}{(2k+1)!
\left|\G(1-\al(k+1))\right|}
$$
\begin{equation}\label{conv unif serie term impares}
\leq \left|\sum_{k=0}^{1}
\frac{-x^{2k+1}}{(2k+1)!\G(1-\al(k+1))}\right|+\sum_{k=2}^{\infty}\frac{|x|^{2k+1}k!}{(2k+1)!\pi}\,.
\end{equation}
 Again, this is an absolutely convergent series in $\bbC$.

  From  (\ref{conv unif serie term pares}) and (\ref{conv unif serie term impares}),
$$ \lim_{\al \nearrow 1}M_{\al/2}(x)$$
$$=\sum_{k=0}^{\infty} \lim_{\al \nearrow
1}\frac{x^{2k}}{(2k)!\G(-\frac{\al}{2}2k+1-\frac{\al}{2})} +
\sum_{k=0}^{\infty}\lim_{\al \nearrow 1}
\frac{-x^{2k+1}}{(2k+1)!\G(1-\al(k+1))}$$
$$
=\sum_{k=0}^{\infty} \frac{x^{2k}}{(2k)!\G(-k+\frac{1}{2})}
=\frac{1}{\sqrt{ \pi}}e^{-\frac{x^2}{4}}.$$
 Moreover, convergence is uniform over compact sets in the
variable $x$.
\endproof 

\begin{theorem}\label{conv a la erf}
  If $x\in \bbR^+_0$ and $\al \in (0,1)$,
$$ \lim_{\al\nearrow 1}\left[1-W\left(-x,-\frac{\al}{2},1\right)\right]
= \mbox{\rm erf\,}\left(\frac{x}{2}\right).$$
\end{theorem}

\proof Observe that
$$ \lim_{\al\nearrow 1}\left[1-W\left(-x,-\frac{\al}{2},1\right)\right]= \lim_{\al\nearrow 1}\int_0^xM_{\al/2}(t)dt$$
and apply Lemma \ref{conv M al/2 cuando al tiende a 1}.
\endproof

\begin{lemma}\label{M pos y decrec}
The Mainardi function $M_{\al/2}(x) $ is a decreasing positive
function if  $\,  0<\al<1$.
\end{lemma} 

\proof Note that  $ M_{\al/2}(x)=W\left(
-x,-\frac{\al}{2},1-\frac{\al}{2} \right)$.

 From \cite{Stankovic} we know that
$$x^{\be -1}W\left( -x^{-\sigma},-\sigma,\beta \right)>0,\, \text{ if }  \quad x>0,   \, \be>0, \, \, 0<\sigma<1, $$
then
\begin{equation}\label{teo stankovic}
W\left( -x^{-\sigma},-\sigma,\beta \right)>0, \, \quad \text{ if }
\,  \quad x>0,   \, \be>0, \, \, 0<\sigma<1.
\end{equation}
In our case, $ \sigma=\frac{\al}{2}\, \in (0,1), \ \be
=1-\frac{\al}{2}>0$, and $g(x)=x^{-\sigma} $ is a one to one
function in $\bbR^+$, so
$$ M_{\al/2}(x)>0 \quad \text{ if } x>0. $$
On the other hand, $
M_{\al/2}(0)=\frac{1}{\G\left(1-\frac{\al}{2}\right)}>0\, , \,
\lim\limits_{x\rightarrow \infty} M_{\al/2}(x)=0$  and
$$ \left(M_{\al/2}(x)\right)'=-W\left(-x,-\frac{\al}{2},1-\al
\right)<0 \,  \text{, because we can apply  (\ref{teo stankovic})
again.} $$ Then the lemma is proved.
\endproof

\begin{corollary}\label{1menosW>0}
If $\,  0<\al<1$ , $1-W\left(-x,-\frac{\al}{2},1\right)$ is a
positive and increasing function in $\bbR^+$.
\end{corollary} 

\proof
It is obvious from $$\left(1-W\left(-x,-\frac{\al}{2},1\right)\right)'=M_{\al/2}(x)>0 \quad
\text{ and }\quad  1-W\left(0,-\frac{\al}{2},1\right)=0.\quad\quad\quad\quad\quad$$
\endproof

\medskip

Is known that (see \cite{FM-The fundamental solution})
\begin{equation}{\label{sol gral para DFE}}
u(x,t)=\int^{\infty}_{-\infty} \frac{t^{-\frac{\al}{2}}}{2\lambda}
M_{\frac{\al}{2}}\left(\left|x-\xi\right|\lambda^{-1}t^{-\frac{\al}{2}}\right) f(\xi)d\xi
\end{equation}
 is a solution for the problem
\begin{equation}{\label{P1}}
\left\{\begin{array}{lll}
          ^C_{0}D^{\al} u(x,t)=\lambda^2\dfrac{\partial^2u }
          {\partial x^2}(x,t) &   -\infty<x<\infty, \,  t>0, \,  0<\al<1, \\
         u(x,0)=f(x) &  -\infty<x<\infty .
      \end{array}\right.\end{equation}
Using this fact, is easy to see that
\begin{equation}{\label{sol DFE 1er cuad}}
 u_1(x,t)=\frac{1}{2\lambda t^{\frac{\al}{2}}}\int^{\infty}_{0}\left[ M_{\frac{\al}{2}}\left(\frac{|x-\xi |}{
\lambda t^{\frac{\al}{2}}}-M_{\frac{\al}{2}}\left(\frac{x+\xi}{
\lambda t^{\frac{\al}{2}}}\right) \right)\right]
 f_0 \, d\xi
\end{equation}
is a solution for the problem
\begin{equation}{\label{P1}}
\left\{\begin{array}{lll}
          ^C_{0}D^{\al} u_1(x,t)=\lambda^2\dfrac{\partial^2u }{\partial x^2}(x,t) &   0<x<\infty, \,  t>0, \,  0<\al<1 \\
         u_1(x,0)=f_0 &   0<x<\infty \\
         u_1(0,t)=0 & t>0.                                  \end{array}\right.\end{equation}

\smallskip
\noindent
Working with (\ref{sol DFE 1er cuad}),
$$u_1(x,t)=\frac{1}{2\lambda t^{\frac{\al}{2}}}\int^{\infty}_{0}\left[ M_{\frac{\al}{2}}\left(\frac{|x-\xi |}{
\lambda t^{\frac{\al}{2}}}-M_{\frac{\al}{2}}\left(\frac{x+\xi}{
\lambda t^{\frac{\al}{2}}}\right) \right)\right]
 f_0 d\xi $$
$$=\frac{f_0}{2}\left[\int^{x}_{0}\frac{1}{\lambda t^{\frac{\al}{2}}}M_{\frac{\al}{2}}\left(\frac{x-\xi }{
\lambda
t^{\frac{\al}{2}}}\right)d\xi+\int^{\infty}_{x}\frac{1}{\lambda
t^{\frac{\al}{2}}}M_{\frac{\al}{2}}\left(\frac{\xi-x }{ \lambda
t^{\frac{\al}{2}}}\right)d\xi\right.$$
$$\left.-\int^{\infty}_{0}\frac{1}{\lambda
t^{\frac{\al}{2}}}M_{\frac{\al}{2}}\left(\frac{x+\xi}{ \lambda
t^{\frac{\al}{2}}}\right)d\xi\right]$$
 $$=\frac{f_0}{2}\left[-W\left(-\frac{x}{\lambda t^{\frac{\al}{2}}},-\frac{\al}{2},1\right)+2 -W\left(-\frac{x}{\lambda
 t^{\frac{\al}{2}}},-\frac{\al}{2},1\right)\right]$$
$$= f_0\left[1-W\left(-\frac{x}{\lambda
t^{\frac{\al}{2}}},-\frac{\al}{2},1\right)\right],
 $$
 and it is easy  to check that
\begin{equation}\label{sol W}
u_2(x,t)=g_0 \, W(-\frac{x}{\la t^{\al/2}},-\frac{\al}{2},1)
\end{equation}
is a solution for the problem
\begin{equation}{\label{P1}}
\left\{\begin{array}{lll}
         ^C_{0}D^{\al} u_2(x,t)=\lambda^2\dfrac{\partial^2u }{\partial x^2}(x,t) &   0<x<\infty, \,  t>0, \,  0<\al<1 \\
         u_2(x,0)=0 &   0<x<\infty \\
         u_2(0,t)=g_0 & t>0.                                  \end{array}\right.\end{equation}

\section*{Acknowledgements}

 This paper has been sponsored by  Project PICTO AUSTRAL
2008 No 73 from {\it Agencia Nacional de Promoci\'{o}n
Cient\'{\i}fica y Tecnol\'{o}gica de la Rep\'{u}blica Argentina}
(ANPCyT) and Project ING349 {\it "Problemas de frontera libre con
ecuaciones diferenciales fraccionarias"}, from Universidad Nacional
de Rosario, Argentina.

We appreciate the valuable suggestions by the anonymous referees
which improved the paper.





 \bigskip \smallskip

 \it

 \noindent
$^1$ Departamento de Matem\'{a}tica - ECEN \\
Facultad de Cs. Exactas, Ingenier\'ia y Agrimensura \\
Universidad Nacional de Rosario \\
Av. Pellegrini 250\
(2000) Rosario, ARGENTINA  \\[4pt]
e-mail: sabrinaroscani@gmail.com
\\[12pt]
 \noindent
$^2$ (Corresp. author) Departamento de Matem\'{a}tica \\
Facultad de Cs. Empresariales\\
Universidad Austral Rosario \\
Paraguay 1950\
(2000) Rosario, ARGENTINA  \\[4pt]
e-mail: edus@fceia.unr.edu.ar



\end{document}